\def\sqr#1#2{{\vcenter{\hrule height.#2pt
      \hbox{\vrule width.#2pt height#1pt \kern#1pt
         \vrule width.#2pt}
       \hrule height.#2pt}}}
\def\eop{\mathchoice\sqr34\sqr34\sqr{2.1}3\sqr{1.5}3}

\def\N{\hbox{\rm I\hskip -0.14em N}}
\def\R{\hbox{\rm I\hskip -0.14em R}}
\def\mod{\hbox{\rm mod}\hskip 0.2em}
\def\th #1 #2. #3\par{\medbreak{\bf#1 #2.
\enspace}{\sl#3}\par\medbreak}
\documentclass{llncs1}
\usepackage{setspace}
\usepackage{makeidx}  % allows for indexgeneration
\usepackage{graphicx}
\begin{document}
\frontmatter          % for the preliminaries
\pagestyle{headings}  % switches on printing of running heads

\title{On the existence of polynomial-time algorithms to the subset sum problem}
\titlerunning{Existence of polynomial-time algorithms}  % abbreviated title (for running head)
\author{Jorma Jormakka}
\authorrunning{Jorma Jormakka}   % abbreviated author list (for running head)
%
%%%% modified list of authors for the TOC (add the affiliations)
\tocauthor{Jorma Jormakka}

\institute{ Contact by: \email{jorma.o.jormakka@gmail.com}}

\maketitle              % typeset the title of the contribution

\begin{abstract}
This paper proves that there does not exist a
polynomial-time algorithm to the the subset sum problem.
As this problem is in $NP$, the result implies that the class $P$
of problems admitting polynomial-time algorithms does not
equal the class $NP$ of problems admitting nondeterministic
polynomial-time algorithms.   
\end{abstract}

\begin {keywordname}
computational complexity, polynomial-time, algorithm, knapsack problem.
\end{keywordname}

\section{Introduction}
Let $\N$ and $\R$ indicate natural and real numbers respectively.

\th Definition 1. A knapsack is a pair of the form 
$(j, (d_1 ,\dots , d_n))$ where
$j,n\in \N$, $j,n>0$ and $d_k\in \N$, $d_k>0$ for $1\le k \le n$.

\vskip 1em
The knapsack problem means the following: given a 
knapsack $(j,(d_1,\dots , d_n))$ determine if there exist 
binary numbers 
$c_k\in \{0,1\}$, $1\le k \le n$, such that
$$j=\sum_{k=1}^nc_kd_k.$$

Let $B,\alpha\in \R$, $B\ge 1$, $\alpha\ge 0$ be fixed numbers. 
An algorithm $A$ is called polynomial-time algorithm 
to the knapsack problem if there exist numbers
$C, \beta\in \R$ that depend on $B$ and $\alpha$ but
not on $n$ such that the following condition is true:
For any sequence of knapsacks of the form
$$((j_n, (d_{1,n},\dots , d_{n,n})))_{n\ge 1}$$ 
satisfying
$$\log_2 j_n<Bn^{\alpha}, \ \log_2 d_{k,n}<Bn^{\alpha}, \ (1\le k\le n), (n\ge 1)\eqno(1.1)$$
the number $N_n$ of elementary operations 
that the algorithm $A$ needs to produce an answer {\it yes} or {\it no} 
to the question if there
exists binary numbers $c_{k,n}\in \{0,1\}$, $1\le k \le n$, such that
$$j_n=\sum_{k=1}^nc_{k,n}d_{k,n}\eqno(1.2)$$
satisfies $N_n<Cn^{\beta}$ for all $n\ge 1$.

The problem that has been described is used in the 
Merkle-Hellman knapsack cryptosystem and today it is commonly known as
the knapsack problem. The name Subset sum problem is used for it 
in \cite {kreh} p. 301, while the name Knapsack problem is reserved
for a more general problem involving selecting objects with 
weights and profits. The name {\it knapsack} is more convenient than 
{\it subset sum} and it is ofen used in this paper.

In the definition of a polynomial-time algorithm for the 
knapsack problem we have included an upper bound on $j_n$ and on each 
$d_{k,n}$, $1\le k\le n$. Such bounds are necessary for the following 
two reasons (i) and (ii).

(i) The number $m$ of bits in the binary representation of $j_n$ satisfies
$m\le \log_2 j_n<m+1$. Thus, if $\log_2 j_n$ grows faster than any polynomial
as a function of $n$ then so does the length of $j_n$ in the binary 
representation. It is necessary to
verify that (1.2) is satisfied. It requires making some operations (like 
compare, copy, read, add, subtract, multiply, divide, modulus) that act on
a representation of $j_n$ on some base number. We may assume that the number
base is 2 as changing a number base does not change the character of 
the algorithm from polynomial-time to non-polynomial-time. 
Any operations that require all bits of $j_n$ 
must require more than a polynomial number of elementary operations 
from any algorithm $A$ if the number of bits in $j_n$ grows faster than 
any polynomial. Similar comments apply to $d_{k,n}$. 

(ii) If $j_n$ has an upper bound independent of $n$, then there exist
a polynomial-time algorithm solving the knapsack problem. 
The Annex gives one such algorithm in Lemma A2.
The algorithm in Lemma A2 calculates an exponentially growing number of 
combinations of $c_k$ in the same polynomial time run. 

Because of (i) and (ii) $j_n$ must grow polynomially with $n$. We can select
$j_n$ as growing linearly as in (1.1). It gives an NP-complete knapsack 
problem.  

{\bf Remark 1.} Lemma A2 in the Annex solves all possible values of 
$j_n<Bn^{\alpha}$ with the same polynomial time 
run of Algorithm A0 because $j_n$
is not used in A0 before checking the final result $b_{n,k}$. 
Let us consider the case when $j_n$ is not limited from above 
by a polynomial of $n$. Lemma A1 runs in polynomial time even 
if the upper bound for $j_n$ grows faster than a polynomial of $n$ but it 
does not produce results that
can tell if there exists a solution for a particular value $j_n$.  
A polynomial-time test, such as taking 
a modulus in (A1), maps the superpolynomial set of possible values 
of $j=\sum_{k=1}^nc_{k,n}d_{k,n}$ into a polynomial number of classes. 
In (A1) the classes are all sums $j$ with the same moduli by $r_n$. 
At least one such a class corresponds to an superpolynomial number
of values $j$. In order to check if any value $j$ in the class equals 
$j_n$ the algorithm should in some way check all of the values $j$ in 
the class, but if the algorithm at the same run checks all values of $j$
then it should in some way loop over a superpolynomial set which
 is not possible for a polynomial time algorithm. 
In general, we can say that a single polynomial time run of an algorithm cannot
solve all values of $j_n$ that are below a superpolynomial upper bound
because the algorithm can only produce a polynomial number of results and
there exist a superpolynomial number of possible values $j_n$.
A polynomial time algorithm that solves the subset sum problem for any value
$j_n$ below a superpolynomial upper bound must limit search and there must
be values $j_n$ that are solved with different runs of the algorithm.  

{\bf Remark 2.} An algorithm is a finite set of rules that at every step
tell what to do next. We can implement an algorithm as a computer program in 
a second generation language on a von Neumann machine and a polynomial time 
algorithm can be implemented in this way so that it requires time and memory 
that grow polynomially with respect to the problem dimension. In the case 
when the smallest upper bound of $j_n$ in Remark 1 grows exponentially
a program in a second generation computer language implementing a 
polynomial time algorithm needs to limit search
by branching instructions, or by acting differently on different data (like in
add, subtract and compare instructions). Thus, we can find values of $j_n$ 
such that the algorithm uses different branches, or acts differently on data,
in solving the subset sum problem.

\section{The inequality (2.6) means non-polynomial time}

It is not possible to select a fixed sequence of 
specific subset sum problems and show that no algorithm can solve this 
specific sequence of problems in polynomial time. This is so because we can 
create an algorithm that treats these specific problems in a particular way 
and can solve that specific  
sequence of problem in a fast way. Instead, we must first select the algorithm
and pose that selected algorithm a sequence of subset sum 
problems that are particularily hard for that specific algorithm. As the 
algorithm can be any possible algorithm, the sequence of problems can only 
be defined by using
some suitable definition of a difficult problem to the selected algorithm and
we cannot give any numerical values for all of the numbers $c_{k,n}$ in
(1.2). The selection will be done by using the following definition of the 
computation time of a subset sum problem.

For convenience, let us select $n$ to be of the form $n=2^{i+2}$ 
for some $i>0$.
This simplifies expressions since it is not necessary to truncate 
numbers to integers.

{\bf Definition 2.}
  
We define a function $f(n)$ that describes (in a certain sense) the worst 
computation time for a selected algorithm.

Let the worst in the median $n$-tuple as be defined as follows. Let 
$$h(d_{1,n},\dots ,d_{n,n},j_n)$$ be the computation time for deciding if
the knapsack $$(j_n,(d_{1,n},\dots ,d_{n,n}))$$ has a solution or not.
Let 
$$Median_{j_n} \ h(d_{1,n},\dots ,d_{n,n},j_n)\eqno(2.1)$$ 
be the median
computation time where $j_n$ ranges over numbers  
$$j_n\in \{ C+1 ,\dots ,2^{n+1}-1\}\eqno(2.2)$$
satisfying the two conditions 
$$j_{n,l}=j_n-C\left\lfloor {j_n\over C}\right\rfloor>2^{{n\over 4}+2}\eqno(2.3)$$
where $C=2^{{n\over 2}+1}$, 
and that there is no solution to the knapsack
$(j_n,(d_{1,n},\dots ,d_{n,n}))$. That is, $j_{n,l}$ are the lower half
bits if $j_n$. 
The values of $j_n$ are computed 
separately in calculation of the median, i.e., 
no partial results from previously computed 
values of $j_n$ are used.

Let
$(d_{1,n},\dots ,d_{n,n})$ range over all knapsack sequences with
$$\lceil \log_2 \sum_{k=1}^n d_{k,n}\rceil=n$$ 
and 
$d_{k,n}\le {2^{n}-1\over n}$.
Because of this requirement at most every second value of $j_n$ in (2.2)
is a solution to the knapsack, i.e., there are $2^n$ combinations of
$(c_{1,n},\dots , c_{n,n})$ mapped to numbers from zero to $2^{n+1}-1$. 
The worst in the median
tuple for $n$ is an $n$-tuple $(d_{1,n},\dots ,d_{n,n})$
(possibly not unique) that maximizes the 
median computation time (2.1). 

Let this maximal median computation time
be denoted by $f(n)$. Thus
$$f(n)=\max_{d_{1,n},\dots ,d_{n,n}} Median_{j_n} \ h(d_{1,n},\dots ,d_{n,n},j_n).\eqno(2.4)$$

We use the median in Definition 2 
instead of the worst case or the worst in the average
case because we need ${n\over 2}$ almost as long computations as the 
worst in (2.6). In the worst and in the worst in the average, a very slow
computation of one value $j_n$ can be the reason for the long computation 
time. By using the median we can find many values $j_n$ giving almost the
median computational time because the distribution of the computational time
for $j_n$ becomes almost normally distributed when $n$ grows due to the
law of large numbers. We include only unsuccessful cases of $j_n$ in the 
computation of the median because this choice implies that a more complicated
knapsack problem (i.e., more cases to check) gives a longer computation time.
If there are more cases to choose, there are more successful cases. Therefore
the time for finding a solution decreases if there are more cases to check.

\th Lemma 1. 
Let $m$ be fixed and $n$ be a power of $m$. 
If $f(n)$ satisfies the inequality
$${n\over m}f\left({n\over m}\right)<f(n)\eqno(2.5)$$
then $f(n)$ does not grow polynomially with $n$.

\proof
Iterating we get
$${n\over m}{n\over m^2}f\left({n\over m^2}\right)<f(n)$$
and iterating up to $k$ yields
$${n^k\over m^{\sum_{i=1}^k i}}f\left({n\over m^k}\right)<f(n)$$
i.e.,
$$e^{ {k\ln n} - {1\over 2}k^2 \ln m -{k\over 2}\ln m} f\left({n\over m^k}\right)<f(n).$$
Setting $k={\ln n\over \ln m}$ (i.e., $1={n\over m^k}$) gives
$$\left( n^{\ln n}\right)^{1\over 2\ln m}n^{-{1\over 2}}f(1)<f(n).$$
If $m$ is any fixed number we see that $f(n)$ satisfying (2.5) is not
bounded by a polynomial function of $n$. $\eop$

\th Lemma 2. 
Let $n$ be a power of $2$. 
If $f(n)=f_1(n)+f_2(n)$ where $f_1(n)$ is a polynomial function of $n$
and $f_2(n)$ satisfies the inequality
$${n\over 2}f_2\left({n\over 2}\right)<f_2(n)\eqno(2.6)$$
then $f(n)$ does not grow polynomially with $n$.

\proof
If $f(n)$ is a polynomial function of $n$ and since $f_1(n)$ is a polynomial
function of $n$ by assumption, it follows that $f_2(n)$ must also be a 
polynomial function of $n$. By Lemma 1, $f_2(n)$ is not a polynomial 
function of $n$, thus neither is $f(n)$. $\eop$

\section{Construction of a special subset sum problem}

In this section we will define a special subset sum problem $K_{1,j_n}$ 
in Definition 3 and show that it can only be solved by solving 
$n_1=n/2$ subknapsacks $(j'_i,(d_{1,n},\dots , d_{n_1,n}))$ 
with different values of $j'_i$.
We will use the denotation $n_1=n/2$ throughout this article
for brevity.  

{\bf Definition 3. Construction of $K_{1,j_n}$.}
We first make a knapsack where the only solutions must satisfy the
condition that exactly one $c_k$ must be $1$ and the others must be zero for
$k=n_1+1$ to $k=n$. 
Let us construct the values 
$d_{k,n}$, $k=n_1+1,\dots ,n$ of $K_{1,j_n}$ for a given $j_n$.
Let $C=2^{{n\over 2}+1}$ and
$$j_{n,h}=C\left\lfloor {j_n\over C}\right\rfloor \ , \  j_{n,l}=j_n-j_{n,h}\eqno(3.1)$$
be the high and low bit parts of $j_n$. Because of (2.2), $j_{n,h}\not=0$. 
Let
$$d_{n_1+k,n}=j_{n,h}+a_k\eqno(3.2)$$
where $0<a_k<\min \{ j_{n,l}, {2^{n_1}-1\over n_1}\}$ 
are distinct integers and there exists no solution
to the knapsack problem for the knapsack 
$$(j'_i,(d_{1,n},\dots , d_{n_1,n}))$$
where 
$$j'_i=j_{n,l}-a_i.\eqno(3.3)$$
Let us also require
that the computation time for $j'_i$ is at least as long as the median
computation time $f(n_1)$ for $(j,(d_{1,n},\dots , d_{n_1,n}))$.
We can select $j'_i$ filling this condition because half of the 
values $j$ 
are above the median. Notice that we compute the median only over values 
$j$ that do not give a solution to the knapsack.
We will also assume that the $j'_i$ are in the set corresponding to
(2.2)-(2.3) for $f(n_1)$, i.e., 
$$j'_i\in \{ C'+1 ,\dots ,2^{n_1+1}-1\}\eqno(3.4)$$
satisfying the condition 
$$j'_i-C'\left\lfloor {j'_i\over C'}\right\rfloor>2^{{n\over 8}+2}\eqno(3.5)$$
where $C'=2^{{n_1\over 2}+1}$.
We may assume so because there are enough values from which 
to choose $j'_i$.

In (3.2) we select the numbers $a_k$ in such a way that the $d_{n_1+k,n}$
satisfy the size condition $d_{n_1+k,n}\le {2^n-1\over n}$. 
Because of the bound (2.3) we have an exponential number of choices for $a_i$.
It is possible to find numbers $j'_i$ such that there is no solution
since only for about half of the values of $j$ there exists a solution
for $(j,(d_{1,n},\dots , d_{n_1,n}))$.
If $j_{n,l}$ is too small and we cannot find values $j'_i$, we take a carry
from $j_{n,h}$ in (3.3) and reselect $a_k$. 
Because of the lower bound on $j$ in (2.2), $j_{n,h}$
is not zero and we can take the carry. Then $j_{n,h}$ is decreased by the 
carry. 

Exactly one $c_k$ must be $1$ and the others must be zero for
$k=n_1+1$ to $k=n$. There cannot be more values $c_k=1$ for $k>n_1$ because 
then the higher bits of $j_n$ are not matched.
The unknown algorithm can try also other combinations but these are the 
only possible combinations and the algorithm must also try them 
(i.e., check these cases in some way unknown to us).  
The sum of the numbers $d_{k,n}$, $k\le {n\over 2}$ is less than
$2^{{n\over 2}+1}-1$. Adding one $c_k$ can give a carry and there may 
not be a solution to the knapsack because the high bits of $j_n$ do
not match but this is not an issue since we do not want solutions. 
We select the $n$-tuple so that there are no solutions to the knapsack
already because the lower bits do not match.

\th Lemma 3.
The algorithm cannot stop to finding a solution
because for every $j_n$ none of the ${n\over 2}$ values of $j'_i$ 
solve the knapsack problem.
Every value $j'_i$ gives at least as long computation as the median
computation time $f(n_1)$.

\proof
We have selected $K_{1,j_n}$ such that 
$(j'_i,(d_{1,n},\dots , d_{n_1,n}))$ has no solution for any $j'_i$.
Thus the algorithm cannot stop because it finds a solution. By construction
the values $j'_i$ give at least as long computation time as the
median for the tuple at $k=1,\dots, n_1$. Since that tuple is the
worst in the median tuple for $n_1$, the computation time 
for each $j'_i$ is at least $f(n_1)$. $\eop$

\th Lemma 4. There is no way to discard any values $j'_i$ without
checking if they solve the subknapsack from $k=1$ to $k=n_1$.
Any case of using the values of $d_{k,n}$ in order to get the result
is considered checking.

\proof
We can select any $a_k$ in such a way that there either exists a 
solution or does not exist. Knowledge from other $c_{i,n}$ ($i\not=n_1+k$)
cannot give any information on how this $a_k$ was selected. 
Thus, the existence
of a solution must be checked using the value $d_{n_1+k,n}$. $\eop$

\th Lemma 5. Several values of $j'_i$ cannot be evaluated on the same run.
The median computation time of $K_{1,j_n}$ is at least
$$f_1(n_1)+n_1f_2(n_1)$$
where $f(n)=f_1(n)+f_2(n)$ is a lower bound for the computation time of one 
$j'_i$ and $f_1(n)$ is a polynomial function of $n$, the shared part of the
computation time of all $j'_i$.

\proof
As explained in Remark 1, a polynomial time algorithm cannot solve all values
of $j'_i$ at the same run because it would require an exponential amount of 
memory. As explained in Remark 2, we can assume that the algorithm is 
implemented in a second generation computer language on a von Neumann machine
and its code has branching instructions, or it acts differently on different 
data in an instruction (like add depends on the data), which has the same 
effect as a branching instruction: for a different $j_n$ there is needed 
a different run. These branching instructions define
a branching tree describing the execution of the algorithm for any input data.
The tree is fixed when the algorithm is selected. At each branching point 
the input data is divided into a finite number of classes. Because this
division is fixed, we can always find two values $j'_i$ which are not executed
by the same polynomial time run. After finding two, we can continue to find 
three values $j'_i$ which all are executed by different polynomial time 
runs of the algorithm. This can be extended to ${n\over 2}$ values $j'_i$: 
we can select $j'_i$ in such a way that no two values $j'_i$ are computed in 
the same run. 
The runs for different values $j'_i$ can have parts that are shared, as long
as the shared parts are computed in polynomial time. 
This is necessarily the case for practical algorithms:
the runs must share at least the beginning of the code 
before branch instructions are reached and this shared part must take only
polynomial time for the algorithm to make any sense. 
The shared part of the computation time can be described by a 
polynomial function $f_1(n)$ and a lower bound for the nonshared computation 
time can be denoted by a function $f_2(n)$.
$\eop$

\section{Proving the inequality (2.6)}
Let the algorithm be chosen. We selected a tuple $K_{1,j_n}$ for a chosen $j_n$
and showed in Lemma 5 
that the computation time for the set of $K_{1,j_n}$ for the single value 
$j_n$ is at least as high as the left hand side of (2.6). 
We have obtained the left side of the inequality (2.6) for an arbitrarily 
chosen algorithm solving the knapsack problem.
However, the set of $K_{1,j_n}$ is a (reasonably) hard problem only for the
chosen value $j_n$. Let us call this $j_n$ with the name $j_{n0}$. 
In the right side of (2.6) 
the number $j_n$ must range over all values and we calculate the median 
computation time over those values of $j_n$ where there is no solution.
In $K_{1,j_{n0}}$ it is very fast to conclude that most values for $j_n$ 
do not have a solution: it is usually enough just to check the bits
of $j_n$ in the most significant half of the number. If they do not match
the most significant bits of $j_{n0}$, then there is no solution.   

We want to change the knapsack problem $K_{1,j_{n0}}$ to another knapsack 
problem $K_2$ (the problem $K_2$ will be defined later in Definition 5)
where $j_n$ can range over all numbers and for many values of $j_n$ there
is no solution and the knapsack problem is difficult. 
The knapsack problem $K_2$ has at most as long median computation time as 
the worst in the median tuple for $n$ because the worst is the worst.

We will do the change in two steps. First we change $K_{1,j_{n0}}$ to
$K_{3,j_{n0}}$ where the bits in the lower half of $j_n$ can vary. In the
second step we change $K_{3,j_{n0}}$ to $K_2$ where also the upper half
bits of $j_n$ can vary.
What we have to show is that the computaton time of the set $K_{1,j_{n0}}$ 
with a single $j_n=j_{n0}$ is not larger than the median computation time
for $K_{3,j_{n0}}$ when $j_n$ can have any lower half bits. 
In $K_{3,j_{n0}}$ only one $d_{j,n}$, the one with $j=n$, has the most 
significant bits of $j_{n0}$.
Therefore $c_{n,n}$ must be one in order to have a possibility of finding
a solution for $j_n$ that has the high bits of $j_{n0}$. We put some numbers
to $d_{j,n}$ for $j=n_1+1,\dots, n-1$. These numbers have zero high bits.
There are more combinations that can give a solution in $K_{3,j_{n0}}$ than
in $K_{1,j_{n0}}$, thus it is easier (and faster) to find a solution, 
provided that there is a solution for a chosen $j_n$. The trick here is that
in the calculation of the median computation time we take only those $j_n$
where there is no solution. Then the fact that there are more possible
combinations only makes it harder to conclude that there is no solution.
We conclude in Lemma 6 that the median computation time for $K_{3,j_{n0}}$ 
when the lower half bits of $j_n$ vary is larger than the computation time
of $K_{1,j_{n0}}$. 

Next we have to show that $K_2$ gives a larger median computation time
when $j_n$ varies over all numbers than $K_{3,j_{n0}}$ when the bits of the
lower half of $j_n$ vary. It is a similar situation here: there are more
combinations in $K_2$ that can give a solution for a given $j_n$, but only
those $j_n$ that give no solution are counted in the median computation time.
Therefore adding complexity makes the median computation time longer. In 
$K_2$ we replace $d_{n,n}$ of $K_{3,j_{n0}}$ by a difficult knapsack problem
in the upper half bits. As this difficult knapsack problem in the upper half
has $n$ numbers $d_{j,n}$ and the bit length of each $d_{j,n}$ is only 
$n/2$, there usually always are solutions to the upper half knapsack problem.
Looking at the upper half knapsack problem does not help in finding 
values $j_n$ that give no solution to the knapsack problem $K_2$. Because
of this, the knapsack problem $K_2$ is not any easier than the knapsack
problem $K_{3,j_{n0}}$. 

Figure 1 shows the main idea. 
\begin{figure}[htbp]
\centering
%\centerline{\includegraphics[scale=0.7,width=5.00in,height=4.00in]{Kuva1a_p_np.jpg}}
\centerline{\includegraphics{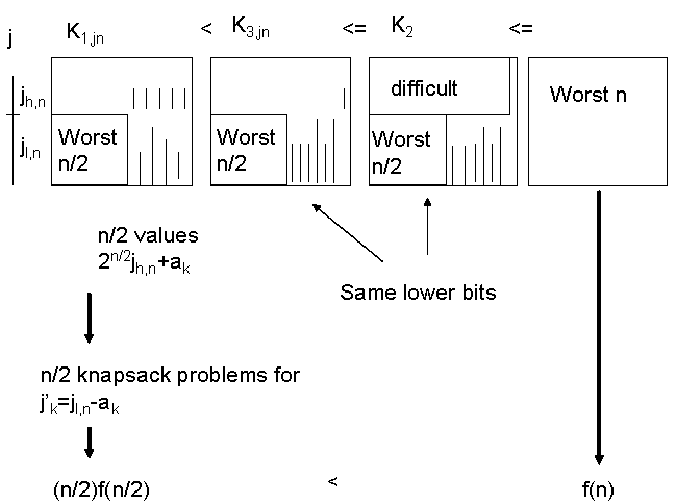}}
\caption{The idea of the proof.}
\label{fig1}
\end{figure}

In Figure 1 the set $K_{1,j_{n0}}$ has the worst
in the median $n_1$-tuple in the left side and the right side has numbers 
from which it is necessary to select exactly one in order to satisfy the
high bits of $j_{n0}$. This yields $n_1$ separate subset sum problems and 
we get the computation time corresponding to the left side of (2.6). 
The set $K_{3,j_{n0}}$ has only one element which has high order bits and it must always be selected in order to satisfy the high bits of $j_n$. Here the bits
of the upper half of $j_n$ are the same as in $j_{n0}$. 
There is the same worst in the median $n_1$-tuple and the remaining $n_1-1$ 
elements can be assigned in any way yielding of the order $n^2$ knapsack 
problems. It is easier to find a solution than in $K_{1,j_{n0}}$, but it
is harder to conclude that there are no solutions.
Lemma 6 shows that the time of solving $K_{1,j_n}$ is not higher than the 
median computation time for $K_{3,j_n}$ for almost any $j_n$ that does not
yield a solution. 

The $n$-tuple $K_2$ has some difficult upper 
half knapsack problem which has to be satistifed with the same values $c_k$
as the lower half knapsack. It is not of any use to check if the upper half
knapsack half has a solution when trying to show that there is no solution 
to the whole knapsack since there almost always are many solutions to the 
upper half knapsack problem. The algorithm must look at all bits.
As finding a solution in $K_2$ requires looking at both the upper and lower
half bits, it should be more difficult to conclude that there are no 
solutions. We will show that at least it is not faster. 
Finally, the inequality from $K_2$ to the 
worst in the median $n$-tuple is obtained directly by the definition 
of what the worst means.

{\bf Definition 4. Construction of $K_{3,j_{n0}}$.}
Let $j_n$ be given and let us define a $n$-tuple $K_{3,j_{n0}}$ 
as an $n$-tuple with elements
$(d_{1,n,3},\dots ,d_{n,n,3})$ 
by specifying the elements
$$d_{k,n,3}=d_{k,n} \hskip 1em (k=1,\dots ,{n\over 2})$$
$$d_{k,n,3}=e_1 \hskip 1em (k={n\over 2}+1,\dots, {3n\over 4})\eqno(4.1)$$
$$d_{k,n,3}=e_2 \hskip 1em (k={3n\over 4}+1,\dots, n-1)$$
$$d_{n,n,3}=j_{n0,h}.$$
We select two nonnegative integers $e_i\le {2^{n_1}-1\over n_1}$, 
$i=1, 2$. The selected $e_1$ and $e_2$ are so small that if $c_n=0$ the higher 
bits of $j_n$ are not matched because there is no carry.
That is, the worst in the median knapsack for $n_1=n/2$ is still in the left
side. The high bits of $j_{n0}$ are in $d_{n,n,3}$. We choose some 
numbers to the elements $d_{k,n,3}$ for $k=n_1+1,\dots, n-1$. 

This $n$-tuple has a simple upper half tuple.
The sum of the numbers $d_{k,n,3}$, $k\le {n\over 2}$ is less than
$2^{{n\over 2}+1}-1$. 
It is always necessary to set $c_n=1$ and this satisfies the upper 
half bits of $j_n$ when $j_n$ ranges over numbers that have the same
upper half bits as $j_{n0}$.

{\bf Definition 5. Construction of $K_2$.}
We will define $K_2$ as an $n$-tuple with elements
$(d_{1,n,2},\dots ,d_{n,n,2})$. Let us remember that the
$n$-tuple $(d_{1,n},\dots ,d_{{n\over 2},n})$ is the worst
in the median tuple for ${n\over 2}$. Let $(d_{0,1},\dots,d_{0,n})$
be an $n$-tuple where each $d_{0,k}\le {2^{n_1}-1\over n_1}$. 
We define 
$$d_{k,n,2}=Cd_{k,n,2}+d_{k,n}\eqno(4.2)$$
for $k=1,\dots , n_1$. 
The numbers $e_1$ and $e_2$ are as in $K_{3,j_{n0}}$ and we define
the elements of $K_2$ for $k=n_1+1$ to $k=n$ as
$$d'_k=Cd_{0,k}+e_1 \hskip 1em (k={n\over 2}+1,\dots, {3n\over 4})$$
$$d'_k=Cd_{0,k}+e_2 \hskip 1em (k={3n\over 4}+1,\dots, n-1)\eqno(4.3)$$
$$d'_n=Cd_{0,n.}$$
Thus, $K_2$ has the same lower half tuple elements as $K_{3,j_n}$ and
in the upper half there is the $n$-tuple $(d_{0,1},\dots,d_{0,n})$.
In this definition we do not specify the $n$-tuple $(d_{0,1},\dots,d_{0,n})$,
but it will be chosen as a sufficietly difficult $n$-tuple.

In $K_{3,j_n}$ our chosen algorithm may fast find a solution and stop
for any $j_n$, but we are only interested at such $j_n$ that give no solution.
The tuple $K_2$ can be split 
into two $n$-tuples: the lower half tuple with 
elements smaller than $C$ and
the upper half tuple that has the higher bit parts. 
In $K_2$ the algorithm
usually does not stop to a solution of the lower half tuple since 
the upper half tuple is usually not satisfied by $c_k$ that satisfy the 
lower half knapsack.

\th Lemma 6. The time for the chosen algorithm to solve $K_{1,j_{n0}}$
is not larger than the median computation time for the algorithm 
for solving $K_{3,j_{n0}}$ when $j_n$ ranges over all values where 
$j_{n,h}=j_{n0,h}$.

\proof
In $K_{3,j_{n0}}$ the indices $k>n_1$ give ${(n+4)n\over 16}$ values of $j$
for a knapsack problem in the indices $k\le n$. Let us name these values $j'_i$
where $i$, $i=1,\dots, {(n+4)n\over 16}$. 

In the indices $k=1,\dots, n_1$ there is the worst in the median $n_1$-tuple. 
The values $j'_i$ that we get are a sample of all possible values $j_{n_1}$ 
for the knapsack problem for this worst in the median $n_1$-tuple. 

Half of all possible values of $j_{n_1}$ yield a longer computation time 
than $f(n_1)$ in the worst in the median knapsack problem for $n_1$
because $f(n_1)$ is the median computation time. If the values
of $j'_i$ that we get are a representative sample of all $j_{n_1}$, then about 
half of the values of $j'_i$ that do not give a solution yield a longer 
computation time than $f(n_1)$. 

We can select $e_1$ and $e_2$ from an exponential set of numbers. Therefore
we can assume that the numbers $j'_i$ are sufficiently well randomly 
distributed over the possible range of the numbers $j_{n_1}$
for the knapsack problem for $n_1=n/2$ and they
are a representative sample of all numbers $j_{n_1}$.

Also, because the numbers $j'_i$ are sufficiently randomly distributed over 
all possible values of $j_{n_1}$ we may assume that about half of the values 
$j'_i$ are on the range (3.4).

There are more values $j'_i$ to check in $K_{3,j_{n0}}$ than the $n/2$ in
$K_{1,j_{n0}}$. If there is no solution for some $j_n$, then it is necessary
to check all $j'_i$ before the algorithm can conclude that there are no 
solutions. Therefore the computation time of the chose algorithm to
solve $K_{1,j_{n0}}$ is not longer than the median computation time for the
algorithm to compute $K_{3,j_{n0}}$ when $j_n$ ranges over all numbers that
have $j_{n,h}=j_{n0,h}$. $\eop$

The median computation time in (2.1) is calculated over the {\it no} 
instances only. Thus, {\it yes} instances are ignored. 
It is sufficient that there are at least 
some {\it no} instances so that (2.1) can be calculated. 
We give an argument that estimates the number of solutions to the knapsack 
problem $(j_n,K_2)$. The argument makes use of averages but it is quite 
sufficient for showing that there are some {\it no} instances for 
computation of (2.1) if the upper bits of $K_2$ are selected in a suitable 
way, indeed a random selection of these bits is likely to yield many 
{\it no} instances.

\th Lemma 7. There are in average $2^{n\over 2}$ solutions 
possible choices of $(c_1,\dots ,c_{n})$ that give the same sum
$\sum_{k=1}^{n}c_kd_{o,k}$.

\proof
The number of combinations of $c_k$ is $2^n$ and the sum  
$\sum_{k=1}^{n}d_{o,k}$ is at most $2^{n\over 2}$.
There are fewer combinations that yield very small or large sums
and most sums are in the middle ranges. $\eop$

\th Lemma 11. We can select the numbers $d_{o,k}$ in such a way that
there are in average about $2^{n\over 4}$ solutions 
possible choices of $(c_1,\dots ,c_{n_1})$ that give the same sum
$\sum_{k=1}^{n_1}c_kd_{o,k}$.

\proof
Most random selections of the numbers $d_{o,k}$ give this result. 
There are fewer combinations that yield very small or large sums
and most sums are in the middle ranges. $\eop$

\th Lemma 12. The lower half tuple in the indices $k=n_1+1,\dots ,n$ 
has only ${n+4\over 4}{n\over 4}$ possible values $j$.

\proof
These numbers are  
$$j=\sum_{k=n_1+1}^nc_k(d_{k,n,2}-Cd_{0,k})=k_1e_1+k_2e_2\eqno(4.4)$$
where $0\le k_1\le {n\over 4}$ and $0\le k_2\le {n\over 4}-1$.
$\eop$

The elements in the worst in the median tuple for $n_1$ satisfy
$d_{k,n}\le {2^{n_1}-1\over n_1}$ 
because we only consider such values of $d_{k,n}$ when finding the worst 
in the median tuple for $n_1$.
Also $e_i\le {2^{n_1}-1\over n_1}$. Thus, there is no carry from the
lower half tuple to the upper half tuple.

\th Lemma 13. It is possible to compute the median (2.1) for $K_2$.

\proof
Let us assume that the values $c_k$ are fixed for the indices 
$k>n_1+1$. This fixes some value $j$ that must be obtained 
from the knapsack in the indices $k=1,\dots ,n_1$ as the subset sum.
By Lemma 12 there are only ${n+4\over 4}{n\over 4}$
possible values $j$.   
The upper half tuple yields about $2^{n\over 4}$ 
possible solutions for a given $j$ in the indices $k=1,\dots ,n_1$ 
by Lemma 11. 
The worst in the median tuple in the lower half tuple has ${n\over 2}$
elements, thus $2^{n\over 2}$ possible numbers can be constructed as sums
$\sum_{k=1}^{n_1}c_kd'_k$ in the lower half tuple.  
The set of the about $2^{n\over 4}$ possible solutions of the 
upper half tuple for a randomly selected $j$
is a small subset of all possible combinations of $c_k$ in the lower
half tuple in the indices $k=1,\dots ,n_1$. 
The probability that any of the 
possible solutions from the upper half tuple is a solution of the 
lower half tuple is only on the range of
${(n+4)n\over 16}2^{-{n\over 4}}$.
The events of selecting the upper half tuple, the lower half tuple, and
the value $j$ can all be considered independent events.
There are only a polynomial number of sums (4.4), thus when $j_n$ is
selected, there are only a polynomial number of possible values for
the lower half of $j$ in $(j,(d'_1,\dots, d'_{n_1}))$.
For a randomly selected $j_n$ there are then only a polynomial number
of $c_k$, $k\le n_1$, that satisfy the lower half bits of $j_n$.
The choice of $c_k$, $k\le n_1$, fixes the upper half of $j$. 
We are left with an upper half knapsack problem for the indices 
$k=n_1+1,\dots ,n$. In this knapsack problem the elements have 
the size about $2^{n_1}$ and there are $n_1$ elements. Thus, for a randomly
selected $j_n$ we expect about one solution.
The solution is constrained by the demand that the 
lower half bits give $j$, i.e., not all combinations are possible.  
We conclude that we get at least some {\it no} instances for 
computation of (2.1) for some choice of $(d_{0,1},\dots ,d_{0,n})$.
$\eop$

\th Lemma 14. The time for the chosen algorithm to solve $K_{3,j_{n0}}$
when $j_n$ ranges over numbers satisfying $j_{n,j}=j_{n0,h}$
is not larger than the median computation time for the algorithm 
for solving $K_2$ when $j_n$ ranges over all values of $j_n$.

\proof
In $K_{3,j_{n0}}$ the upper bits are easily satisfied by selecting 
$c_{n,n}=1$. 
In order to find a solution to the subset sum problem for $K_2$
the algorithm must find a common solution to two knapsacks, i.e., both
the upper bits and the lower bits knapsacks in $K_2$
must be solved with the same numbers $(c_1,\dots ,c_n)$. 
We may choose any difficult knapsack $(d_{0,1},\dots ,d_{0,n})$
to the upper bits of $K_2$. 

The algorithm cannot conclude that there are no solutions to the whole 
knapsack problem because there are no solutions to the upper half knapsack
problem. This is so since there almost always are many solutions to the 
upper half knapsack problem for any value of $j$: the upper half 
knapsack problem has $n$ elements of the bit length at most $n/2$. This means
that there are $2^n$ possible combinations of $c_k$ and they are mapped to
$2^{n/2}$ different numbers $j$. Each number $j$ is likely to come from 
many combinations of $c_k$ since in average $2^{n/2}$ combinations give the 
same $j$. 

It is also not possible to the algorithm to check that none of the
solutions to the upper half knapsack problem give a solution to the lower
half knapsack problem. This is so because there are exponentially many 
(i.e., $2^{n/2}$) solutions to the upper half knapsack problem. They cannot be
checked in a polynomial time. 

Because of these two reasons the median computation
time of $K_{3,j_{n0}}$ when $j_n$ ranges over all $j_n$ that has the 
same high bits as $j_{n0}$ cannot be higher than the median computation time 
for $K_2$ where $j_n$ ranges over all numbers. In the computation of the median
time we only take cases of $j_n$ where there is no solution and a more 
complicated $n$-tuple must give a longer time for concluding that there
are no solutions. $\eop$

\th Lemma 15. The inequality (2.6) holds for the chosen algorithm.

\proof 
By Lemma 6 the median computation time for $K_{3,j_{n0}}$ 
when the median is taken over the set of $j_n$ having $j_{n,h}=j_{n0,h}$
is at least as high as the time to solve $K_{1,j_{n0}}$.
By Lemma 13 we can calculate the median of computation times over cases
when there is no solution for $K_2$.
By Lemma 14 the median computation time for $K_2$ when $j_n$ ranges over 
all values is not smaller than the median computation time for 
$K_{3,j_{n0}}$ when the median
is computed over the set $j_n$ where $j_{n,h}=j_{n0,h}$.
As $K_2$ is a fixed $n$-tuple it follows from the definition of the worst
in the median tuple that $K_2$ has at most as long median computation
time as the worst in the median tuple for $n$, i.e., $f(n)$.
Thus the inequality (2.6) holds. $\eop$

\bigskip

\th Theorem 1. Let an algorithm for the knapsack problem
be selected. There exist
numbers $B,\alpha\in \R$, $B\ge 1$, $\alpha\ge 0$ and
a sequence
$$((j_n,(d_{1,n},\dots , d_{n,n})))_{n\ge 1}$$
of knapsacks satisfying
$$\log_2 j_n<Bn^{\alpha}, \ \log_2 d_{k,n}<Bn^{\alpha}, (1\le k\le n), \ (n\ge 1)$$
such that the 
algorithm cannot determine in polynomial time if there exist
binary numbers $c_{k,n}$, $1\le k\le n$, satisfying
$$j_n=\sum_{k=1}^nc_{k,n}d_{k,n}.$$

\proof
The idea of this proof is to compare the computation time of the worst 
(in some sense) knapsack of size
$n$ to the computation time of (in the same sense) worst knapsack of 
${n \over 2}$.
The  computation time was defined in (2.4) and denoted by $f(n)$.
By Lemma 15 the inequality (2.6) holds for an arbitrary chosen algorithm. 
By Lemma 2 the arbitrarily chosen algorithm is not a polynomial time 
algorithm. 
$\eop$

\th Theorem 2. {\bf P} does not equal {\bf NP}.

\proof
The knapsack problem is well known to be in {\bf NP}. $\eop$

\section{Annex}

\th Lemma A1. Let $B\ge 1$, $\alpha\ge 0$ and $\gamma\ge 0$ be selected. 
Let $r_n>0$ and $j_n$ be integers satisfying 
$$r_n<n^{\gamma},\hskip 2em \log_2 j_n<Bn^{\alpha}\hskip 3em (n\ge1).$$
There exist numbers $C, \beta\in \R$,$C\ge 1$, $\beta\ge 0$ and
an algorithm that given any sequence of knapsacks
$$((j_n,(d_{1,n},\dots , d_{n,n})))_{n\ge 1}$$
can determine for each $n$ if there exist
binary numbers $c_{k,n}$, $1\le k\le n$, such that
$$j_n\equiv \sum_{k=1}^nc_{k,n} d_{k,n}\hskip 1em (\mod r_n).\eqno(A1)$$
The number $N_n$ of elementary operations needed by the algorithm 
satisfies $N_n<Cn^{\beta}$ for every $n>1$.

\proof
The bound on the logarithm of $j_n$ guarantees that modular arithmetic
operations on $d_{k,n}$ can be made in polynomial time since we can assume
that $d_{k,n}\le j_n$.
We can find the numbers $c_{k,n}$ by computing numbers $s_{k,j,n}$
from the recursion equations for $k$
$$s_{k,j,n}= s_{k-1,j,n}+s_{k-1,(j-d_{k,n})(\mod r_n), n}\eqno(A2)$$
$$s_{0,j,n}=\delta_{j=0},$$
where the index $j$ ranges from $0$ to $r_n-1$ and is calculated modulo $r_n$.
The index $n$ is fixed and only indicates that the numbers are for the 
$n^{th}$ knapsack.
Here $\delta_x$ is an indicator function: $\delta_x=1$ if the statement
$x$ ( i.e., $j$ equals $0$ in (A2) ) is true and $\delta_x=0$ if $x$ is false.
Let
$$G_{k,n}(x)=\sum_{j=0}^{r_n-1}s_{k,j,n}x^j,$$
where $|x|<1$. From (A2) follows
$$\sum_{j=0}^{r_n-1}s_{k,j,n}x^j= \sum_{j=0}^{r_n-1}s_{k-1,j,n}x^j+\sum_{j=0}^{r_n-1}s_{k-1,(j-d_{k,n})(\mod r_n),n}x^j.$$  
Changing summation to $j'=j-d_{k,n}$ yields
$$G_{k,n}(x)= G_{k-1,n}(x)+\sum_{j'=-d_{k,n}}^{r_n-1-d_{k,n}}s_{k-1,j'(\mod r_n),n}x^{j'+d_{k,n}}.$$
Changing the order of summation of $j'$ shows that
$$G_{k,n}(x)= G_{k-1,n}(x)+x^{d_{k,n}}\sum_{j'=0}^{r_n-1}s_{k-1,j',n}x^{j'}.\eqno(A3)$$
Simplifying (A3) gives
$$G_{k,n}(x)= G_{k-1,n}(x)+x^{d_{k,n}}G_{k-1,n}(x).$$
As $G_{0,n}(x)=s_{0,0,n}=1$, we get
$$G_{n,n}(x)= \prod_{k=1}^n(1+x^{d_{k,n}}).$$
Expanding the product shows that 
$s_{k,j,n}\not=0$ if and only if there exist binary numbers $c_m$,
$c_m\in \{0,1\}$, $1\le m\le n$, satisfying
$$j\equiv \sum_{m=1}^nc_md_{m,n}\hskip 2em (\mod r_n).$$
For $j=j_n$ and $k=n$ we get the knapsack problem. This means that we
can solve the knapsack problem by computing all $s_{k,j,n}$ form (A2). 
We do not actually need the numbers $s_{k,j,n}$ but only
the information if $s_{k,j,n}\not=0$. Therefore we will not compute
the terms $s_{k,j,n}$ directly but   
calculate binary numbers $b_{j,k}\in \{0,1\}$ by Algorithm A0 below.
The number $b_{k,j}$ calculated by A0 is zero if and only if
the number $s_{k,j,n}=0$ is zero.
\vskip 1em 
Algorithm A0:

\hskip 1em Loop from $k=0$ to $k=n$ with the step $k:=k+1$ do $\{$

\hskip 2em Loop from $j=0$ to $j=r_n-1$ with the step $j:=j+1$ do 

\hskip 3em $b_{j,k}:=0$

\hskip 1em $\}$

\hskip 1em $b_{0,0}:=1$

\hskip 1em Loop from $k=1$ to $k=n$ with the step $k:=k+1$ do $\{$

\hskip 2em $M:=\min \{r_n-1, \sum_{m=1}^kd_{m,n}\}$ 

\hskip 2em Loop from $j=0$ to $j=M$ with the step $j:=j+1$ do $\{$

\hskip 3em If ($b_{k-1,j}=0$ and $b_{k-1,(j-d_{k,n})(\mod r_n)}=0$) do $b_{j,k}:=0$

\hskip 3em else do $b_{j,k}:=1$

\hskip 2em $\}$

\hskip 1em $\}$

\hskip 1em If $b_{n,j_n}=1$ do $result:=TRUE$ else do $result:=FALSE$

\vskip 1em

Algorithm A0 loops from $k=0$ to $k=n$ and from $j=0$ to 
$j=r_n-1<n^{\gamma}$. Thus A0 needs 
a polynomial number of elementary operations as a function of $n$ in order to
give the result $TRUE$ or $FALSE$ to the existence of a solution to (A1).
$\eop$

\vskip 1em

\th Lemma A2. Let $B,\alpha\in \R$, $B\ge 1$, $\alpha\ge 0$ 
be fixed.
There exist numbers $C, \beta\in \R$, $C\ge 1$, $\beta\ge 0$ and
an algorithm that for any sequence
$$((j_n,(d_{1,n},\dots , d_{n,n})))_{n\ge 1}$$
of knapsacks satisfying
$$j_n\le Bn^{\alpha},\hskip 3em d_{k,n}\le j_n\hskip 3em(1\le k\le n),$$
can determine if there exist
binary numbers $c_{k,n}$, $1\le k\le n$, such that
$$j_n=\sum_{k=1}^nc_{k,n}d_{k,n}.$$  
The number $N_n$ of elementary operations needed by the algorithm 
satisfies $N_n<Cn^{\beta}$ for every $n>1$.

\proof
The result follows directly from Lemma A1 by selecting 
$r_n=\sum_{k=1}^nd_{k,n}\le nj_n$. 
$\eop$

\end{document}